\documentclass[conference]{IEEEconf}
\IEEEoverridecommandlockouts    % to override locked commands

\usepackage{multirow}
\usepackage{lipsum}
\usepackage{amsmath}
\usepackage{amssymb}
\usepackage[ruled,vlined]{algorithm2e}
\usepackage{graphicx}
\graphicspath{{./Figures/}}
\title{\LARGE \bf Hybrid System Identification of Electric Freight Transition Dynamics via SINDy and Neural ODEs}

 \author{
 	\parbox{\textwidth}{%
 		\centering
 		Alessandro Carli$^{1}$, 
        Maret Clement$^{2}$,
        Jonas Mårtensson$^{1}$, Anna Pernestål$^{3}$, Matthieu Barreau$^{1}$%
 	}%
 	\thanks{$^{1}$ Department of Decision and Control Systems, Digital Futures, KTH Royal Institute of Technology, Stockholm, Sweden {\tt\small acarli@kth.se, jonas1@kth.se, barreau@kth.se}}%
    \thanks{$^{2}$ Department of Mechanical Engineering, Eindhoven University of Technology, the Netherlands {\tt\small maretc@kth.se}}%
 	\thanks{$^{3}$ Integrated Transport Research Lab (ITRL), Stockholm, Sweden {\tt\small pernestal@kth.se}}%
 }
\begin{document}

	\maketitle
	\thispagestyle{empty}
	\pagestyle{empty}
	%%%%%%%%%%%%%%%%%%%%%%%%%%%%%%%%%%%%%%%%%%%%%%%%%%%%%%%%%%%%%%%%%%
\begin{abstract}
The transition to electric heavy-duty freight is constrained by a strong interdependence between vehicle adoption and charging infrastructure deployment. While system dynamics models are well-suited to simulate these socio-technical feedback loops, their reliance on conditional logic and heuristic rules makes the underlying dynamics difficult to interpret mathematically. This limits the direct application of formal systems and control theory.
This paper proposes a solution to that problem through a general hybrid system identification framework that extracts a continuous-time analytical surrogate from simulations of complex system dynamics. The approach combines the sparse identification of nonlinear dynamics algorithm with neural ordinary differential equations to form a grey-box model. The first component identifies the dominant interpretable dynamics under causal coupling constraints, while the neural residual captures unmodelled nonlinearities.
The framework is trained using multiple shooting over a 40-year horizon. The resulting model reproduces the training trajectories with a normalized root-mean-square error below 4\%, while maintaining reliable predictive accuracy when evaluated on unseen initial conditions.%
%By extracting explicit governing dynamics, the proposed hybrid framework provides a mathematically tractable foundation for future systems analysis and optimization.
\end{abstract}
	
	%%%%%%%%%%%%%%%%%%%%%%%%%%%%%%%%%%%%%%%%%%%%%%%%%%%%%%%%%%%%%%%%%%

\section{Introduction}
\label{sec:introduction}

Road transportation remains a primary driver of global climate change, contributing approximately 17.8\% of total global energy-related CO\textsubscript{2} emissions \cite{Liu2022}. Within this sector, heavy-duty trucks account for roughly 21\% of emissions \cite{Vaccaro2024}, making the electrification of freight transport a critical strategy for decarbonization \cite{Yuan2021}. While the passenger vehicle market has experienced rapid growth, increasing its market share from less than 5\% in 2020 to over 20\% in 2024, the electric heavy truck sector has not yet seen a comparable acceleration. Current sales shares remain low, reaching only 4.4\% in China and 2.2\% in Europe \cite{InternationalEnergyAgency2025}.

A major complication in this transition is the circular dependency between vehicle adoption and infrastructure development. The high power requirements of heavy trucks necessitate substantial infrastructure investment, yet charging providers are hesitant to deploy capital without a guaranteed fleet size to ensure profitability. Conversely, freight operators are unlikely to electrify their fleets without a mature and reliable charging network. To analyse these complex adoption dynamics and evaluate the potential impact of policy interventions, researchers frequently rely on System Dynamics (SD) models. This study utilizes an underlying Vensim SD framework developed in \cite{raoofi2025electric}, which captures the nonlinear feedback loops, time delays, and multi-stakeholder interactions inherent in the transition to electric freight.

However, while this existing SD framework is well-suited for scenario exploration, it acts as a mathematical gray box. Its dynamics are encoded through discrete lookup tables, complex internal procedures, and conditional logic, rather than explicit governing equations. Consequently, policy evaluation in this environment relies on manual parameter tuning and repeated forward simulations, which are computationally inefficient and limit systematic optimization. More fundamentally, even when SD frameworks export numerical routines, the resulting models remain algorithmically opaque. The lack of a closed-form symbolic representation prevents the analysis using formal tools from systems and control theory, such as nonlinear stability analysis, sensitivity bounds, or optimal control.
This disconnect between simulation-based insight and mathematical analysis remains a well-known challenge in large-scale, interconnected system models \cite{Poort2025, Schilders2008}.

%To better analyze these types of complex systems and overcome the opacity of traditional simulations, recent advances in data-driven system identification have emerged. Sparse regression methods, in particular, have demonstrated strong performance in recovering nonlinear dynamics from complex energy and infrastructure systems \cite{Nandakumar2023}. In parallel, the emergence of scientific machine learning has enabled the systematic integration of physics-based structure with data-driven models through the concept of universal differential equations \cite{Rackauckas2020}. Furthermore, combining sparse regression with neural representations has been shown to improve robustness when modeling noisy or partially observed systems \cite{Chen2021}. Building on these advancements, this paper proposes a hybrid, data-driven system identification pipeline to bridge the gap between simulation-based logistics planning and formal control theory. By fusing Sparse Identification of Nonlinear Dynamics (SINDy) with Neural Ordinary Differential Equations (Neural ODEs) \cite{Chen2018}, we extract a grey-box analytical surrogate from the Vensim model. This surrogate is more accurate than its white-box equivalent, more expressive, and more interpretable than its pure black-box counterpart. Indeed, it retains the core dynamic properties of the original simulation while remaining mathematically tractable for future control and stability analysis. 
To better analyse these complex systems and overcome the opacity of traditional simulations, recent advances in data-driven system identification have emerged. Sparse regression methods, in particular, have shown strong performance in recovering nonlinear dynamics from complex energy and infrastructure systems \cite{Nandakumar2023}. In parallel, scientific machine learning has enabled the integration of physics-based structure with data-driven models through universal differential equations \cite{Rackauckas2020}. Furthermore, combining sparse regression with neural representations improves robustness when modelling noisy or partially observed systems \cite{Chen2021}. Building on these advances, this paper proposes a hybrid data-driven system identification pipeline to bridge the gap between simulation-based logistics planning and formal control theory. By combining Sparse Identification of Nonlinear Dynamics (SINDy) with Neural Ordinary Differential Equations (Neural ODEs) \cite{Chen2018}, we extract a grey-box analytical surrogate from the Vensim model. This surrogate is more accurate than its white-box equivalent, while remaining more interpretable than a pure black-box model. It preserves the core dynamic properties of the original simulation while remaining mathematically tractable for future control and stability analysis.

The paper is structured as follows: Section \ref{sec:Background} presents the model and summarizes the limitations of standard system identification approaches. Section \ref{sec:Method} introduces the proposed hybrid SINDy and Neural ODE pipeline. Section \ref{sec:Results} reports and analyses the results, focusing on training fidelity and validation on previously unseen datasets. Section \ref{sec:conclusion} concludes with a discussion of implications for policy optimization and stability analysis.

\begin{figure*}
    \centering
    \includegraphics[width=0.75\linewidth]{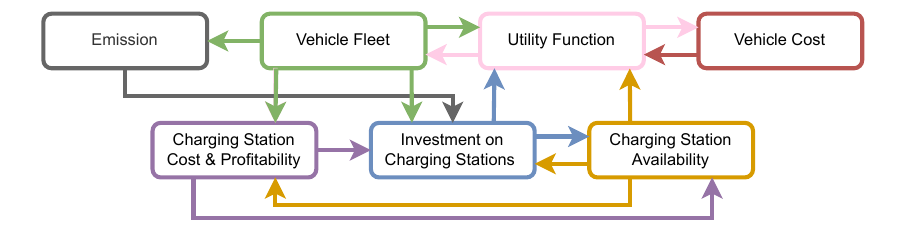}
    \caption{Visual representation of the seven sectors in Vensim and their connections (from \cite{raoofi2025electric}).}
    \label{fig:Vensim_submodules_overview}
\end{figure*}

%%%%%%%%%%%%%%%%%%%%%%%%%%%%%%%%%%%%%%%%%%%%%%%%%%%%%%%%%%%%%%%%%%
\section{Background and Problem Formulation}
\label{sec:Background}

To resolve the mathematical opacity of the existing simulation, the system must be translated into a set of explicit ordinary differential equations (ODEs). Traditional system identification often struggles with the high dimensionality and nonlinearity typical of socio-technical systems, while purely deep-learning-based alternatives lack interpretability and struggle to generalize outside their training distribution \cite{Rackauckas2020}. 

Therefore, it is more common to use SD to model this kind of system. In the case discussed in this paper, the transition to electric freight, the SD model has been developed with Vensim through seven interconnected sectors, as illustrated in Fig.~\ref{fig:Vensim_submodules_overview}. Each of these sectors is interpretable, but their interconnection is creating the system complexity.

\begin{figure}[h!]
    \centering\includegraphics[width=0.95\linewidth]{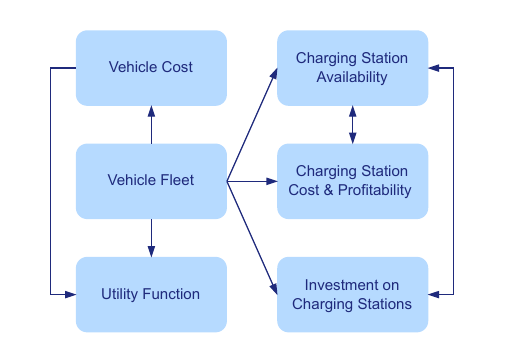}
    \caption{Schematic of the coupled sectors within the model. Arrows indicate the principal causal feedback loops. This structural knowledge informs the sparsity constraints utilized in the proposed identification method.}
    \label{fig:coupled_sectors}
\end{figure}

Within this framework, our analysis specifically focuses on the primary state variables, represented as stocks in the SD model, across six sectors, visualized in Fig.~\ref{fig:coupled_sectors}, with inputs $u_i$ and outputs $y_i$ for $i \in \{1, \dots, 6\}$. We obtain trajectory data of inflows and outflows to each sector for $N_p$ different initial conditions and parameters denoted $\{p_j\}_{j\in\{1, \dots, N_p\}}$, leading to the following dataset:
\begin{equation}
    \label{eq:dataset}
    \mathcal{D}_{N_p}= \left\{ \left( k, \left[ \begin{matrix} u_1^j(k) \\ \vdots \\ u_6^j(k) \end{matrix} \right], \left[ \begin{matrix} y_1^j(k) \\ \vdots \\ y_6^j(k) \end{matrix} \right] \right) \right\}_{\substack{k\in\{2017,\dots,2060\}\\ j\in\{1, \dots, N_p\}}}.
\end{equation}
The problem addressed in this paper is to propose a methodology for automatically deriving an analytical formulation of the interconnection between these six sectors based on the knowledge of $\mathcal{D}_{N_p}$.

Before moving on to the method itself, we will derive the main properties of the SD models that will be leveraged later on. An initial observation indicates that the dynamics of each of these six sectors are based on stock variables that represent accumulations within the system. Their values evolve over time, driven by the net difference between their respective inflows and outflows. Because the temporal evolution of these accumulations is fundamentally determined by rates of change, the behaviour of the simulation can be mathematically expressed in the standard second-order continuous-time dynamical formulation $\dot{x}_i = f_i(x_i, u_i)$ where $x_i = [ x_{i,1} \ x_{i,2}]$, $f_i$ represents the dynamics of sector $i$, and $u_i$ are the inputs to this sector, and $y_i = h_i(x_i)$ are the outputs. In this work, we will consider that $h$ is the identity function and that we directly measure the stock variables.

This inherent mathematical structure provides the basis for applying the SINDy algorithm. While this method is effective at recovering standard continuous dynamics, the SD model incorporates non-analytic structures that sparse polynomials struggle to capture: fixed time delays in charging infrastructure decay ($t - \tau$), discrete lookup tables for investor sensitivity (e.g., piecewise profitability indices), and hard nonlinearities in government-to-private investment ratios. Furthermore, nested conditional logic (\texttt{IF-THEN-ELSE}) creates discrete switching based on temporal thresholds (2030, 2045) and emission goals.

To address these limitations, we propose a sequential hybrid architecture that derives a reduced-order, grey-box surrogate of the electric truck adoption system. This framework fuses Sparse Identification of Nonlinear Dynamics \cite{Brunton2016SINDy}, used to recover interpretable, dominant physical structures, with Neural ODEs \cite{Chen2018}, which serves as a residual network. By embedding the interpretable SINDy equations within a differentiable Neural ODE framework, we extract the explicit governing dynamics of the adoption process where possible, while relying on deep learning to approximate the remaining complex residuals. The resulting model is a computationally efficient and formally tractable surrogate.
%suitable for future stability analysis and policy optimization.

\section{Method}
\label{sec:Method}

To accurately recover the governing dynamics of the electric truck adoption system, we propose a fully integration-based, hybrid identification pipeline. The methodology extracts the dominant physical structure via an integration-based variant of SINDy (SINDy-RK4) and subsequently captures unmodelled residuals via a Neural ODE.

\subsection{Structural Recovery (SINDy-RK4)}

As established in Section \ref{sec:Background}, the evolution of the system's stock variables is driven by continuous rate flows. Consequently, the electric truck adoption process is initially modelled as a continuous-time nonlinear dynamical system with control inputs:
\begin{equation}
    \left\{ 
    \begin{array}{l}
        \dot{x}(t) = f(x(t), u(t)), \\
        x(t_0) = x_0,
    \end{array} 
    \right.
    \label{eq:ODE}
\end{equation}
where $x(t) = \left[ \begin{matrix} x_{1,1} & x_{1,2} & \cdots & x_{6,1} & x_{6,2} \end{matrix} \right] \in \mathbb{R}^{12}$ represents the state vector of the whole system, $u(t) \in \mathbb{R}^m$ denotes the input variables, and $f$ is the flow function. We denote by $t \mapsto \pi_i(x_0, t)$ the trajectory of sector $i$ with \eqref{eq:ODE} starting at $x_0$.

Following the SINDy framework \cite{Brunton2016SINDy} and its formal extension for actuated systems (SINDYc) \cite{Kaiser2018, Kaiser2016}, we assume that the underlying vector field $f$ is sparse in the space of possible mathematical functions. This assumption is justified by the structure of the underlying simulation: stock variables are updated by specific, localized flow rates rather than a dense web of dependencies, and the interconnections between sectors is sparse. Enforcing this sparsity is critical, since it prevents overfitting the training data and ensures the resulting model remains interpretable by isolating only the dominant, active physical mechanisms.

To approximate the unknown function ${f}$, an expansive dictionary of candidate nonlinear functions, ${\Theta}({x}, {u})$, is constructed. This physical feature library is designed to capture the relevant governing dynamics and may contain constant terms, linear combinations of the coupled state and control variables, as well as higher-order nonlinear interactions depending on the assumed complexity of the system:
\begin{equation}
    {\Theta}({x}, {u}) = \begin{bmatrix}
        1 & {x} & {u} & {x}^2 & {x}{u} & \dots
    \end{bmatrix} \in \mathbb{R}^{1 \times p},
    \label{eq:library}
\end{equation}
where $x^2$ and $xu$ denote the vectors of element-wise quadratic and cross-interaction terms, respectively, and $p$ represents the total number of candidate features generated. By expanding or restricting this library, the framework can flexibly accommodate various degrees of polynomial nonlinearity, enabling, for instance, Taylor expansions of more complex nonlinear functions.

The identification problem is thus posed as finding a sparse combination of these basis functions that accurately reconstructs the system dynamics:
\begin{equation}
    \dot{x}(t) \approx {\Theta}({x}, {u}){\Xi},
    \label{eq:sindy_approx}
\end{equation}
where ${\Xi} \in \mathbb{R}^{p \times n}$ is the sparse coefficient matrix. Each column of ${\Xi}$ corresponds to a specific state variable, and the non-zero entries within that column dictate exactly which terms from the feature library actively govern its time evolution.

Standard SINDy relies on numerical differentiation to approximate $\dot{{x}}(t)$. However, numerical differentiation acts as a high-pass filter, amplifying any measurement noise or abrupt state transitions characteristic of discrete simulations. To overcome this limitation, we utilize an established integration-based sparse regression formulation, leveraging a 4th-order Runge-Kutta (RK4) scheme as proposed by \cite{Goyal2022}. By embedding this formulation within a differentiable programming framework \cite{Rackauckas2020}, the sparse coefficients ${\Xi}$ are treated as trainable parameters rather than variables in a static algebraic system. 
Optimizing over the integrated trajectory effectively smooths stochastic noise and forces the identified physics to remain stable over long horizons by penalizing the compound growth of local errors.

The dynamics are discovered by minimizing the mean squared error over discrete simulation steps between the measured data and the integrated trajectory generated by the RK4 solver:
\begin{multline}
    \mathcal{L}_{SINDy} = \underbrace{\frac{1}{S} \sum_{k=1}^{S} \frac{1}{N_p} \sum_{j=1}^{N_p} \frac{1}{6} \sum_{i=1}^6 \| \pi_i(p_j, k) - y_i(k) \|^2}_{\text{data loss}} \\ 
    + \lambda \underbrace{\|{\Xi}\|_1}_{\text{sparsity loss}},
    \label{eq:loss}
\end{multline}
where $S$ is the number of sequence steps in the evaluated temporal batch, $\pi_i(p,\cdot)$ is the state trajectory of sector $i$ dynamically generated by the RK4 solver solving \eqref{eq:ODE} with initial condition $p$, and $y$ is the measurement from the collected dataset $\mathcal{D}_{N_p}$ in \eqref{eq:dataset}. The parameter $\lambda$ controls the strength of the $L_1$ regularization. By penalizing the absolute sum of the coefficients, the $L_1$ norm forces the weights of non-essential functions to zero, which promotes model parsimony.

%The variable $\lambda$ is the regularization parameter enforcing sparsity on the coefficient matrix.
 
Training differential equations over long temporal horizons via backpropagation methods is highly susceptible to gradient explosion \cite{Chen2018}. To stabilize the learning of ${\Xi}$ during this phase, we employ a discrete multiple-shooting strategy paired with a curriculum learning schedule. The continuous training data is sliced into randomized temporal batches of length $\tau$. To prevent early divergence, the simulation horizon $\tau$ is linearly increased as training progresses across epochs $e$:
\begin{equation}
    \tau(e) = \min \left( T_{max}, \ \tau_{min} + \alpha \cdot e \right).
    \label{eq:curriculum}
\end{equation}
This formulation forces the optimizer to resolve local dynamics before generalizing to global trajectory tracking \cite{Bengio2009}.
Furthermore, while $L_1$ regularization promotes small coefficients, it rarely yields the exact zeros required for structural equation discovery. To enforce strict sparsity, we employ an iterative hard-thresholding scheme. Periodically during training, any coefficient $|\Xi_{ij}|$ falling below a structural threshold is forced to zero. Once a term is truncated, it is permanently removed from the optimization space by zeroing its corresponding gradient, ensuring the identified model structure remains parsimonious throughout the remaining epochs.

This approach provides a robust sparse relaxation \cite{Zheng2019} that prevents pruned features from re-entering the system during subsequent optimization steps, by forcing the optimizer to minimize trajectory error using only the most relevant physical terms.

\subsection{Structural Constraints: Coupling Matrix}

Directly applying SINDy to a high-dimensional system results in an exponential growth of candidate terms. Recent hybrid architectures demonstrate that embedding physical priors into sparse regression significantly improves robustness and reduces computational cost \cite{Wang2025}. 
To prevent overfitting within the expansive polynomial library ${\Theta}$, a coupling matrix $\mathcal{C} \in \{0, 1\}^{n \times n}$ is introduced as a strong inductive bias. The matrix $\mathcal{C}$ encodes the allowable causal interactions between sectors, derived from the established system dynamics structure. If sector $j$ does not influence sector $i$ ($\mathcal{C}_{ij} = 0$), the corresponding parameters in ${\Xi}$ are never initialised. This reduces the optimization problem from a global search to a set of physically bounded, lower-dimensional problems.

In this study, the coupling matrix $\mathcal{C}$ is configured to focus on the primary interactions between the six selected sectors, resulting in 10 active couplings among the 30 possible cross-state couplings. While the source simulation contains a broader range of variables, we limit our identification to the stock variables that define the system's behaviour. Consequently, connections associated with excluded variables are not represented in the current matrix. 

%we focus on six sectors, yielding 30 possible cross-state couplings. Using a sparse configuration matrix, only 10 interactions are identified as active, corresponding to an approximate 67\% reduction in modelled couplings. Diagonal terms capture intrinsic state dynamics, while sparsity is enforced on off-diagonal coupling terms. Encoding inter-sector couplings in this manner enables policy updates without requiring modifications to the underlying model or computational pipeline.

\subsection{Trajectory Correction via Residual Neural ODEs}
While SINDy-RK4 effectively recovers the dominant physical structure, it is constrained by the fixed basis functions in ${\Theta}$. To capture unmodelled dynamics, such as complex saturation limits or hidden empirical behaviour, a residual neural network is integrated.

The full hybrid system is formulated as an additive correction:
\begin{equation}
    \dot{x}(t) = \underbrace{{\Theta}({x}, {u}){\Xi}}_{\text{Fixed Physics}} + \underbrace{f_{\text{res}}({x}, {u}; {\phi}).}_{\text{Neural Residual}}
    \label{eq:hybrid_dynamics}
\end{equation}
This is trained in a phased approach. In this subsequent phase, the sparse physical coefficients ${\Xi}$ discovered are frozen. The neural network weights ${\phi}$ are optimized to capture the residual error between the fixed physical model and the ground truth.

The residual component $f_{\text{Res}}$ is implemented as a feedforward neural network. It takes a 24-dimensional input, formed by concatenating the state variables and exogenous inputs. The architecture consists of two hidden layers with 32 neurons each, utilizing hyperbolic tangent activation functions, followed by a 12-dimensional linear output layer.

To ensure the model prioritizes the physical explanation, the weights and biases of this final output layer are explicitly initialised to zero. The neural network effectively only corrects the residuals that the interpretable SINDy model cannot explain, preventing the black-box component from dismantling the known physics.

\subsection{Joint Optimization}
Building upon the two-stage identification, we introduce an additional joint refinement phase, the parameter-freezing constraints are lifted to enable holistic refinement. Both the sparse physical coefficients ${\Xi}$ and the neural weights ${\phi}$ are optimized simultaneously. This joint fine-tuning serves as a final calibration step, ensuring strict consistency between the identified physics and the neural correction, ultimately yielding superior stability and predictive accuracy on validation datasets.

\section{Results}
\label{sec:Results}
    
To assess the fidelity and interpretability of the identified hybrid model, we evaluate performance across two dimensions: quantitative simulation accuracy over long-term horizons, and the structural balance between the physical and neural components.

\subsection{Evaluation Metrics}

\subsubsection{Accuracy}
To assess the fidelity and robustness of the identified hybrid model, we evaluated the simulation accuracy over the full time domain. The primary metric, Normalized Root Mean Square Error (NRMSE), quantifies the deviation of the learned trajectory from the ground truth relative to each state's natural variance. The metric is defined as:

\begin{equation}
    \text{NRMSE} = \frac{\sqrt{\frac{1}{S} \sum_{k=1}^{S} \| \hat{x}(t_k) - x_{true}(t_k) \|^2}}{\sigma_{true}},
\end{equation}
where $\sigma_{true}$ represents the standard deviation of the ground truth state trajectory. A low NRMSE serves as the primary indicator of the pipeline's predictive accuracy, confirming that the hybrid model closely tracks the trajectories across the full simulation.
Additionally, we use the coefficient of determination to evaluate the correlation between predicted and observed trends. 

\subsubsection{Hybrid Contribution Ratio ($\rho$)}:
To ensure the structural integrity, we adopt the principle of model parsimony \cite{Brunton2016SINDy}, verifying that the neural component does not overwhelm the interpretable physical terms. We quantify this through the hybrid contribution ratio $\rho$. This metric %effectively a measure of relative residual magnitude \cite{Rackauckas2020}, 
ensures that the neural network acts as a localized corrector and is defined as the magnitude of the neural correction relative to the physical term:
\begin{equation}
    \rho = \frac{\|f_{\text{res}}({x}, {u})\|_2}{\|{\Theta}({x}, {u}){\Xi}\|_2}.
\end{equation}

This metric serves as a proxy for interpretability:
\begin{itemize}
    \item Sparse correction ($\rho \ll 1$): The neural network acts as a minor adjustment to a dominant physical model. This is the desired ``grey-box'' behaviour.
    \item Dominance ($\rho \ge 1$): The neural network dominates the dynamics, indicating that the sparse polynomial basis failed to capture the system structure.
\end{itemize}

\subsection{Training Performance Analysis}

\begin{figure}[ht]
    \centering
    \includegraphics[width=\linewidth]{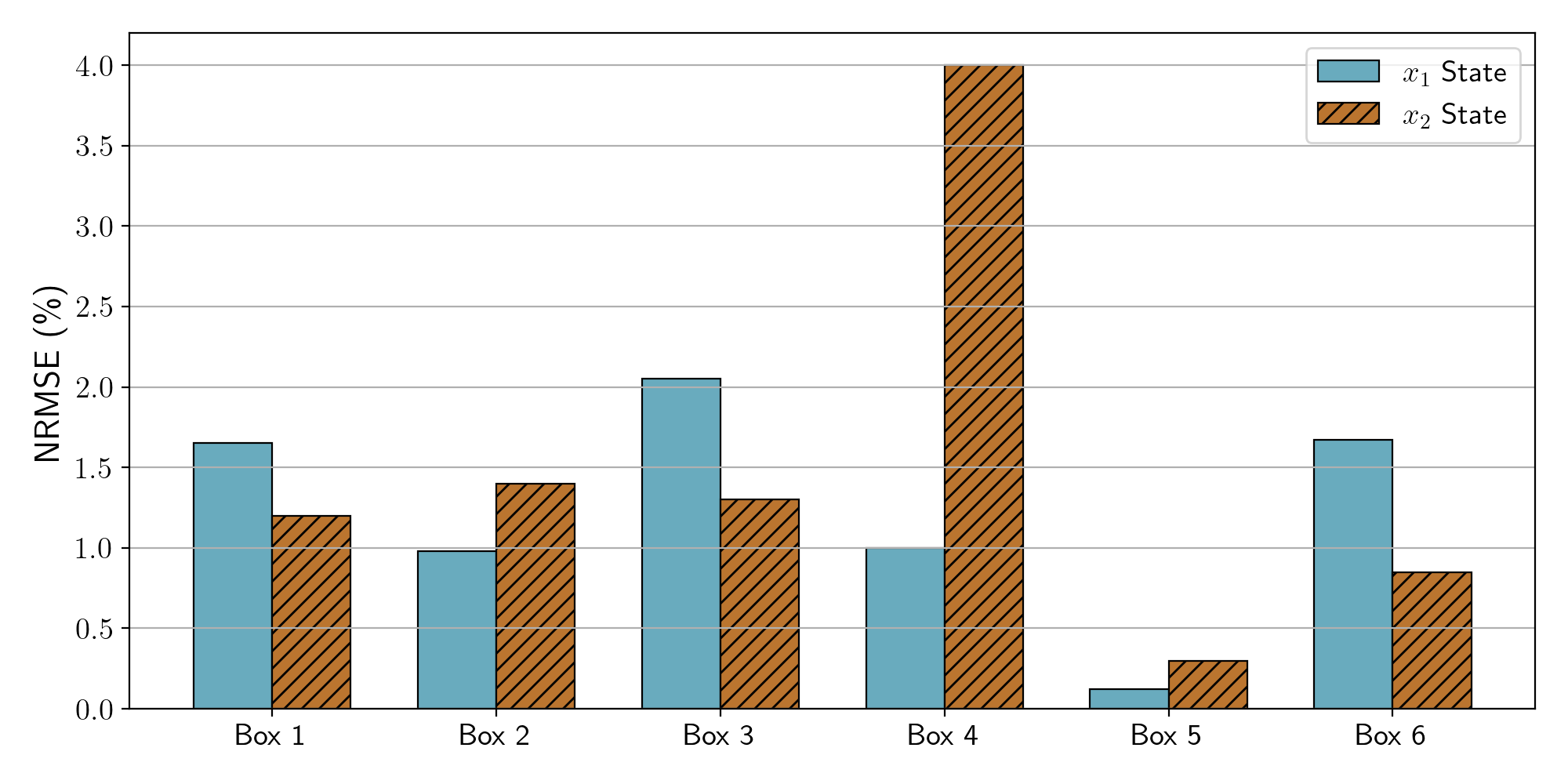}
    \caption{Normalized Root Mean Square Error by sector. The blue bars ($x_1$) and brown dashed bars ($x_2$) represent the primary and secondary states, respectively.}
    \label{fig:nrmse_by_box}
\end{figure}

Fig.~\ref{fig:nrmse_by_box} presents the error decomposition by sector and state variable ($x_1, x_2$), averaged across the seven socio-economic scenarios used during training. 
The aggregate performance indicates that the hybrid model successfully stabilizes the dynamics, with NRMSE values generally bounded below $4\%$. 
These results indicate that the hybrid framework can reliably track coupled socio-technical dynamics, although approximation accuracy naturally varies across sectors.

The simulation accuracy across the coupled system provides critical insight into both the capabilities and the structural limitations of polynomial-based identification. The error profile reveals distinct behavioural clusters that directly correlate with the mathematical formulations of the underlying SD model:

\begin{itemize}
    \item Inactive and strictly continuous dynamics: The states in Box 5 exhibit exceptional accuracy, with an NRMSE below 0.5\%. However, this high fidelity is largely an artifact of the system boundaries; for example, the first state ($x_1$) functions as an inactive accumulator that remains at zero throughout the simulation. Such constant and strictly smooth trajectories are trivially resolved by the sparse polynomial basis identified in the first phase.
    
    \item Non-analytic discontinuities: The highest deviations are localized in state $x_1$ in Box 3 (2.1\% error) and the second state of Box 4 (4\% error). These errors are not random artifacts, but directly reflect the non-smooth formulations inherent in the Vensim architecture. The dynamics in these sectors are governed by explicit conditional logic (\texttt{IF-THEN-ELSE}), bounding operators (\texttt{MAX}), and sharp saturation constraints.
    Polynomial bases, such as ${\Theta}({x}, {u})$, are inherently continuous and globally smooth. They struggle to approximate transitions and flattened saturation regions created by \texttt{MAX} functions without introducing oscillations. At these discrete logical switching points in the vector field $\dot{{x}}$, the SINDy component effectively smooths over the discontinuity. This forces the neural residual to bear the primary modelling burden, and small integration errors inevitably accumulate over the horizon, resulting in the observed drift.
    
    \item Moderate nonlinearity: The remaining sectors exhibit intermediate errors ranging from 1.0\% to 2.0\%. In these regions, the identified polynomial captures the macroscopic trends in the trajectory. The neural residual is not forced to compensate for hard switches; rather, it serves as a minor corrective mechanism to address local discretization and truncation errors.
\end{itemize}

Ultimately, this heterogeneity highlights the utility of the hybrid architecture. While standard sparse regression struggles to capture the discrete logical switches common in socio-technical simulations, the neural ODE effectively absorbs these non-analytic residuals without corrupting the interpretable physical baseline.

To verify that the hybrid model remains physically grounded, we monitored the hybrid contribution ratio $\rho$ throughout the phased training process. During training, $\rho$ was initialised near zero ($0.0036$) due to the zero-weighted final layer of the residual network. As the curriculum learning schedule extended the simulation horizon from 8 to 40 steps, the ratio stabilized at approximately 0.14. This indicates that the neural component contributes less than 12.5\% of the total vector field magnitude. This result confirms that the sparse polynomial basis identified captures the dominant dynamics, while the neural network functions strictly as a localized corrective mechanism for unmodelled nonlinearities.

\subsection{Validation on Unseen Scenarios}

\begin{figure}[ht]
    \centering
    \includegraphics[width=\linewidth]{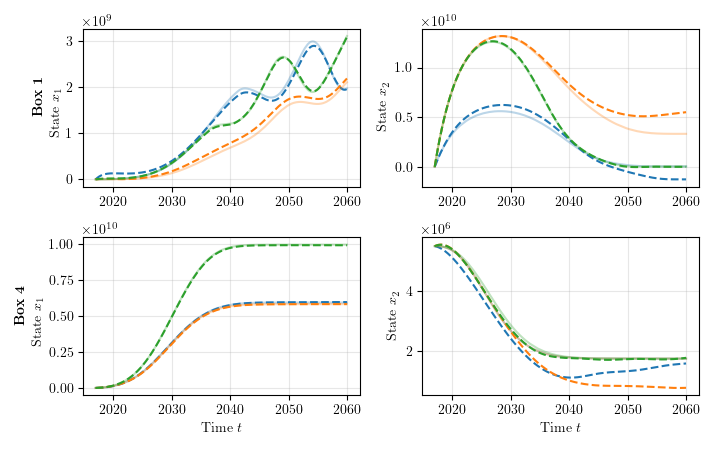}
    \caption{Temporal evolution of selected state variables, comparing the model predictions against ground truth across the three datasets.}
    \label{fig:validation_ts}
    \vspace*{-0.3cm}
\end{figure}

To evaluate the generalization capabilities of the identified model, cross-validation was performed on three distinct datasets representing alternative initialization scenarios. Fig.~\ref{fig:validation_ts} compares the predicted trajectories (dashed lines) against the ground truth (solid lines) over the simulation horizon.

\begin{figure}[ht]
    \centering
    \includegraphics[width=\linewidth]{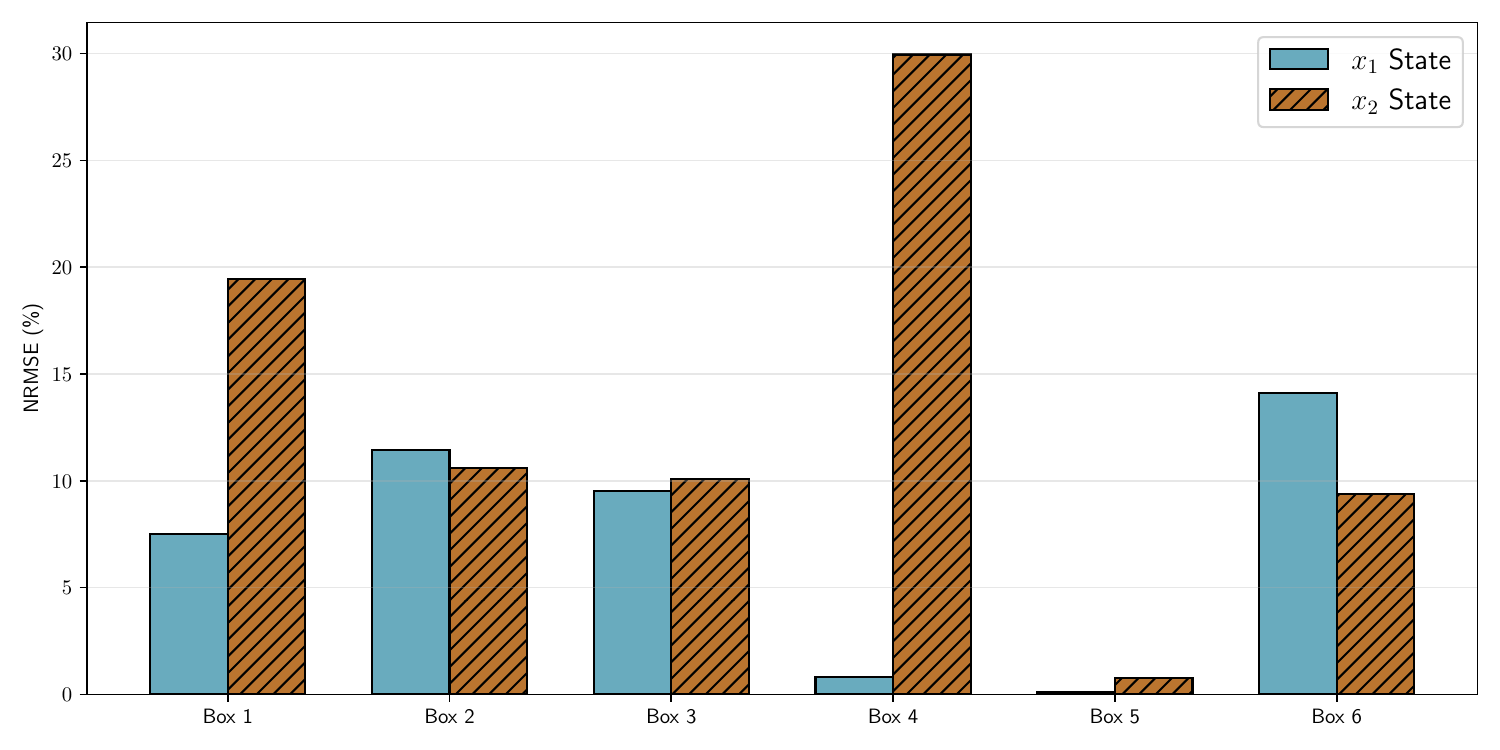}
    \caption{Normalized Root Mean Square Error by sector for validation datasets. The blue bars ($x_1$) and brown dashed bars ($x_2$) represent the primary and secondary states, respectively.}
    \label{fig:validation_rmse}
    \vspace*{-0.3cm}
\end{figure}

Fig.~\ref{fig:validation_rmse} presents the performance across the validation datasets, measuring the scale-invariant accuracy of the predictions through NRMSE. 

The validation results highlight two specific numerical phenomena that reflect the structural properties of the hybrid framework:
\begin{itemize}
    \item 
State $x_1$ in Box 5 yields a negative $R^2$ score alongside a nominal NRMSE of 0.13\%. This is a standard statistical artifact occurring when the truth signal is constant, resulting in near-zero variance. Because the metric is defined as $R^2 = 1 - \frac{SS_{\text{res}}}{SS_{tot}}$, a total sum of squares ($SS_{tot}$) approaching zero causes the ratio to mathematically diverge. The low NRMSE confirms that the model correctly maintains the inactive state without generating inauthentic dynamics.

\item
Consistent with the training observations, Box 4 ($x_2$) remains the most difficult sector to approximate, yielding an NRMSE of 30.0\% and an $R^2$ of 0.88. The time-series plots demonstrate that while the hybrid model captures the initial trajectory decay, it struggles to track the asymptotic tail. This deviation reinforces the idea that the state is governed by a discrete, logical capacity switch.
\end{itemize}

Excluding these structural limits, 9 of the 12 states maintain an $R^2$ above 0.97. This suggests that the hybrid architecture successfully captures the system's primary dynamics and generalizes to unseen initialisation conditions.

%\subsection{Robustness Analysis}

%To confirm that the identified hybrid model learned a generalized, stable vector field rather than merely overfitting the training data, its robustness was evaluated through Monte Carlo simulations using the learned dynamics.

%To check if the model remains stable under noise, additive Gaussian noise (scaled by 0.05) was injected into both the normalized initial conditions (${x}_0$) and the control trajectories (${u}(t)$) prior to integration. Despite the susceptibility of long-horizon neural differential equations to error amplification, the resulting simulated trajectories remained strictly bounded within a 95\% confidence interval ($\pm 2\sigma$) over the simulation. The absence of mathematically divergent simulations confirms that the joint optimization scheme successfully embedded stabilizing physical priors into the network, allowing it to resist noise propagation.

\section{Conclusion}
\label{sec:conclusion}
This study presented a grey-box identification pipeline for extracting a continuous-time surrogate from a system dynamics model of electric truck adoption. By fusing sparse symbolic regression with neural differential equations, we successfully isolated the interpretable physical backbone of the electric truck transition while delegating the simulation's inherent logical discontinuities, such as threshold government interventions and infrastructure delays, to a neural residual. 
While the framework demonstrated high predictive fidelity across diverse scenarios, the identification process highlighted a fundamental trade-off: the model accurately captures broad adoption trends but faces inherent limitations when approximating constrained, asymptotic logic.

The resulting representation provides a formally tractable foundation for the next phase of this research. By moving beyond the black-box nature of system dynamics, this surrogate enables the direct calculation of system sensitivities and the deployment of control methods. Future work will leverage this to perform formal stability analysis on fleet electrification trajectories and to derive optimal policy interventions. Ultimately, it enables a more robust assessment of the feedback loops that define the transition to sustainable freight.
%%%%%%%%%%%%%%%%%%%%%%%%%%%%%%%%%%%%%%%%%%%%%%%%%%%%%%%%%%%%%%%%%%
\section*{ACKNOWLEDGMENTS}
% Replace with acknowledgments or remove if none
%The authors would like to thank [Funding Agency/Person] for support.
This work is partially supported by the Wallenberg AI, Autonomous Systems and Software Program (WASP) funded by the Knut and Alice Wallenberg Foundation, and ITRL (Integrated Transport Research Lab) at KTH.

%%%%%%%%%%%%%%%%%%%%%%%%%%%%%%%%%%%%%%%%%%%%%%%%%%%%%%%%%%%%%%%%%%
%\addtolength{\textheight}{-12cm}
%\vspace{10mm}
\bibliographystyle{IEEEtran}
% Your .bib file here
\bibliography{root} 
	
\end{document}